\documentclass[12pt,leqno]{article}
\usepackage{amsmath, amsthm, amsfonts, amssymb, color}
\usepackage{mathrsfs}
\theoremstyle{definition}

\newcommand{\scr}[1]{\mathscr #1}
\definecolor{wco}{rgb}{0.5,0.2,0.3}

\numberwithin{equation}{section} \theoremstyle{remark}

\newcommand{\ua}{\uparrow}

\title{{\bf Derivative  Formulae and Poincar\'e Inequality for   Kohn-Laplacian Type Semigroups }\footnote{Supported in
 part by NNSFC(11131003) and the Laboratory of Mathematical and  Complex Systems.}
}
\author{
{\bf Feng-Yu Wang$^{1),2)}$ }\\
\footnotesize{$^{1)}$School of Mathematical Sciences,
Beijing Normal University, Beijing 100875, China}\\
  \footnotesize{$^{2)}$Department of Mathematics,
Swansea University, Singleton Park, SA2 8PP, UK}\\
\footnotesize{Email: wangfy@bnu.edu.cn; F.Y.Wang@swansea.ac.uk} }
\begin{document}
\def\R{\mathbb R}  \def\ff{\frac} \def\ss{\sqrt} \def\B{\mathbf
B}
\def\N{\mathbb N} \def\kk{\kappa} \def\m{{\bf m}}
\def\dd{\delta} \def\DD{\Delta} \def\vv{\varepsilon} \def\rr{\rho}
\def\<{\langle} \def\>{\rangle} \def\GG{\Gamma} \def\gg{\gamma}
  \def\nn{\nabla} \def\pp{\partial} \def\EE{\scr E}
\def\d{\text{\rm{d}}} \def\bb{\beta} \def\aa{\alpha} \def\D{\scr D}
  \def\si{\sigma} \def\ess{\text{\rm{ess}}}
\def\beg{\begin} \def\beq{\begin{equation}}  \def\F{\scr F}
\def\Ric{\text{\rm{Ric}}} \def\Hess{\text{\rm{Hess}}}
\def\e{\text{\rm{e}}} \def\ua{\underline a} \def\OO{\Omega}  \def\oo{\omega}
 \def\tt{\tilde} \def\Ric{\text{\rm{Ric}}}
\def\cut{\text{\rm{cut}}} \def\P{\mathbb P}
\def\C{\scr C}     \def\E{\mathbb E}
\def\Z{\mathbb Z} \def\DDD{\mathbf D}
  \def\Q{\mathbb Q}  \def\LL{\Lambda}
  \def\B{\scr B}    \def\ll{\lambda}
\def\vp{\varphi}\def\H{\mathbb H}\def\ee{\mathbf e}

\maketitle
\begin{abstract} As a generalization to  the heat semigroup on the Heisenberg group, the diffusion semigroup generated by the subelliptic operator $L:=\ff 1 2 \sum_{i=1}^m X_i^2$ on $\R^{m+d}:= \R^m\times\R^d$ is investigated, where
$$X_i(x,y)= \sum_{k=1}^m \si_{ki} \pp_{x_k} + \sum_{l=1}^d (A_l
x)_i\pp_{y_l},\ \ (x,y)\in\R^{m+d}, 1\le i\le m$$ for $\si$ an invertible $m\times m$-matrix and $\{A_l\}_{1\le l\le d}$ some $m\times m$-matrices such that the H\"ormander condition holds. We first establish Bismut-type and Driver-type derivative formulae with applications on gradient estimates and the coupling/Liouville properties, which are new even for the heat semigroup on the Heisenberg group; then extend some recent results derived  for the heat semigroup on the Heisenberg group.
\end{abstract} \noindent

 AMS subject Classification:\ 60J75, 60J45.   \\
\noindent
 Keywords: Kohn-Laplacian type operator, derivative formula,   Poincar\'e inequality, reverse Poincar\'e inequality.
 \vskip 2cm

\section{Introduction}

In recent years, the heat semigroup generated by the Kohn-Laplacian
on  the Heisenberg group regularity   has been intensively
investigated, see \cite{BBBC, DM, Li} for derivative
estimates and applications, and see \cite{BB1, BB2, BG, BGM} (where a more
general model was considered) for the generalized curvature
conditions and applications.

The first purpose of this paper is to establish Bismut's derivative
formula \cite{B} and Driver's integration by parts formula \cite{D}
for the semigroup generalized by a class of Kohn-Laplacian type
operators. These two formulae are crucial for stochastic analysis of
diffusion processes and are not explicitly known   even for the heat semigroup on the Heisenberg group.
 Our second aim is to  extend some known results derived recently for the heat semigroup on the Heisenberg group to a more general framework
of Kohn-Laplacian type
operators. 

Let us first recall the Kohn-Laplacian on the three-dimensional
Heisenberg group. Consider the following two vector fields on
$\R^3$:
$$X_1(x)= \pp_{x_1} -\ff{x_2}2 \pp_{x_3},\ \ X_2(x)= \pp_{x_2} +\ff {x_1} 2
\pp_{x_3},\ \ x=(x_1,x_2,x_3)\in\R^3.$$ Then $\DD_K:= X_1^2+X_2^2$ is called the
Kohn-Laplacian. It is crucial in the study of this operator that
$[X_1,X_2]=\pp_{x_3}, [X_i,\pp_{x_3}]=0 (i=1,2)$ and $X_1,X_2$ are left-invariant under the
group action $$(x_1,x_2,x_3)\bullet (x_1',x_2',x_3') =
\Big(x_1+x_1', x_2+x_2', x_3+x_3' +\ff 1 2 (x_1x_2'-x_2x_1')\Big).$$

To do stochastic analysis with this operator, let us introduce the
associated stochastic differential equation for $(X(t),Y(t))\in
\R^2\times \R:$
$$\beg{cases} \d X(t)=   \,\d B(t),\\
\d Y(t) = \<A X(t), \d B(t)\>,\end{cases}$$ where $B(t)$ is the
$2$-dimensional Brownian motion and $A=   \Big(\beg{matrix} 0 &-1\\
1 &0\end{matrix}\Big).$ Then $(X(t),Y(t))$ is the diffusion process
generated by $\ff 1 2\DD_K$, and the associated transition semigroup is known as the heat semigroup on the Heisenberg group.

In this paper we consider the following
natural extension of this equation for
$(X(t),Y(t))\in \R^m\times \R^d=:\R^{m+d}$ $(m\ge 2,d\ge 1)$: \beq\label{E1}\beg{cases} \d X(t)= \si \,\d B(t),\\
\d Y_l(t) = \<A_l X(t), \d B(t)\>,\ \ 1\le l\le
d,\end{cases}\end{equation} where $B(t)$ is the $m$-dimensional
Brownian motion on a complete probability space $(\OO,\F,\P)$ with the natural filtration $\{\F_t\}_{t\ge 0}$, $\si$ is an invertible $m\times m$-matrix, and $(A_l)_{1\le l\le d}$ are $m\times m$-matrices.
Let
$$X_i(x,y)= \sum_{k=1}^m \si_{ki} \pp_{x_k} + \sum_{l=1}^d (A_l
x)_i\pp_{y_l},\ \ (x,y)=(x_1,\cdots, x_m, y_1,\cdots, y_d)\in\R^{m+d}, 1\le i\le m.$$ Then the solution to (\ref{E1}) is
the diffusion process generated by $$ L:=\ff 1 2 \sum_{i=1}^m
X_i^2.$$ Obviously, for any $1\le i,j\le m$ and $ 1\le l\le d$, we
have  $[X_i, \pp_{y_l}]=0$   and
$$[X_i,X_j]= \sum_{l=1}^d
\big\{(A_l\si)_{ji}-(A_l\si)_{ij}\big\}\pp_{y_l}=\sum_{l=1}^d (G_l)_{ji} \pp_{y_l},$$ where $G_l:= A_l\si-\si^*A_l^*.$ Then the H\"ormander
condition holds (thus, $L$ is subelliptic) if and only if
\beq\label{Ho}  {\rm The}\
\{m(m-1)\} \times d\text{-matrix}\  (M_{(i,j),l})_{1\le i,j\le m; 1\le
l\le d}\ {\rm has\  rank}\  d,\end{equation}  where
$$M_{(i,j),l}:= (G_l)_{ij},\ \ 1\le i<j\le m, 1\le
l\le d;$$ or equivalently,
\beq\label{Ho'} \sum_{i,j=1}^m \Big|\sum_{l=1}^d (G_l)_{ij} a_l\Big|^2 \ge \ll |a|^2,\ \ a=(a_l)_{1\le l\le d}\in \R^d\end{equation} holds for some constant $\ll>0.$

A simple example for (\ref{Ho}) or (\ref{Ho'}) to hold is that
$d=m-1, \si=I_{m\times m}$ and \beq\label{A}(A_l)_{ij} =\beg{cases} \aa_l, &\text{if}\ i=1, j=l+1,\\
\bb_l, &\text{if}\ i=l+1, j=1, \\
0,  &\text{otherwise}\end{cases}\end{equation}  for $\aa_l\ne \bb_l, 1\le l\le d.$

Moreover, let $\R^{m+d}$ be equipped with the group action
$$(x,y)\bullet (x',y') =
(x+x', y+y' +\<(\si^*)^{-1}A_\cdot x, x'\>),\ \ (x,y), (x',y')\in\R^{m+d}, $$ where $\<(\si^*)^{-1}A_\cdot x, x'\>:= (\<(\si^*)^{-1}A_l x, x'\>)_{1\le l\le d}\in\R^d.$   Then $(0,0)$  is the unique unit element,   and the inverse element of  $(x,y)\in\R^{m+d}$  is
$$(x,y)^{-1}:= (-x, \<\si^{-1}A_\cdot x,x\> -y).$$ It is easy to see that $\{X_i\}_{1\le i\le m}$ are left-invariant vector fields under the group structure. Indeed, for any $f\in C^1(\R^{m+d})$ and $(u,v)\in\R^{m+d}$, letting $f_{(u,v)}(z)= f((u,v)\bullet z), z\in\R^{m+d}$, we have
\beg{equation*}\beg{split} X_i f_{(u,v)} (0,0) &= \sum_{k=1}^m \si_{ki} \Big\{\pp_{x_k} f_{(u,v)}\Big\}(0,0)=\sum_{k=1}^m \si_{ki} \Big\{\pp_{x_k} f +\sum_{l=1}^d ((\si^*)^{-1} A_l u)_k  \pp_{y_l} f\Big\}(u,v) \\
&= \Big\{\sum_{k=1}^m \si_{ki}   \pp_{x_k} f + \sum_{l=1}^d (A_l x)_i \pp_{y_l} f\Big\}(u,v) = (X_i f)(u,v),\ \ 1\le i\le m. \end{split}\end{equation*} It is also easy to see that the Lebesgue measure $\mu$ is invariant under the group action.

We will investigate the Markov semigroup $(P_t)_{t\ge 0}$ for the solution to equation (\ref{E1}):
$$P_t f(x,y):= \E f(X^x(t), Y^{(x,y)}),\ \ (x,y)\in\R^{m+d},\ t\ge 0, f\in \B_b(\R^{m+d}),$$ where $(X^x(t), Y^{(x,y)}(t))_{t\ge 0}$ is the solution to the equation with initial data $(x,y).$ Since div$X_i=0, 1\le i\le m,$ $P_t$ is symmetric in $L^2(\mu)$.

In Section 2 we investigate  Bismut/Driver-type  derivative formulae for $P_t$ and applications.  In Section 3 and Section 4 we modify the argument in \cite{BBBC} to derive explicit Poincar\'e and reverse Poincar\'e inequalities for $P_t$. As we emphasized in Introduction that explicit derivative formulae are new even for the heat semigroup on the Heisenberg group. Moreover, although functional and Harnack inequalities derived in \cite{BGM} using generalized curvature conditions apply to our present framework,   results derived therein do not cover our Poincar\'e inequality and   explicit inverse Poincar\'e inequality.

 \section{Derivative formulae}

Recall that  $G_l:= A_l \si - \si^* A_l^* (1\le l\le d)$ are skew-symmetric, i.e. $G_l^*=-G_l$. In this section we assume
 \paragraph{(A1)} \ $G_l\ne 0$ for all $1\le l\le d$, and there exists a constant $\theta\in [0,1)$ such that
 $$\theta \sum_{l=1}^d a_l^2 |G_l u|^2 \ge   \sum_{l,k: 1\le l\ne k\le d} |a_la_k \<G_l^*G_k u, u\>|,\ \ u\in\R^m, a\in \R^d.$$

 It is easy to see that {\bf (A1)} implies the H\"ormander condition. Indeed, letting $u={\bf e}_i$, we obtain from {\bf (A1)} that
$$\theta \sum_{l=1}^d (G_l^*G_l)_{ii} \ge \sum_{1\le l\ne k\le d} \big|a_la_k(G_l^*G_k)_{ii}\big|,\ \ 1\le i\le m, a\in\R^d.$$ Therefore, for any $a=(a_l)_{1\le l\le d}\in\R^d$,  we have 
\beg{equation*}\beg{split} &\sum_{1\le i,j\le m} \Big|\sum_{l=1}^d M_{(i,j),l}a_l\Big|^2 =\sum_{k,l=1}^d\sum_{i,j=1}^m(G_l)_{ij}(G_k)_{ij}a_la_k\\
&= \sum_{k,l=1}^d \text{Tr}(G_k^*G_l) a_k a_l \ge (1-\theta) \sum_{l=1}^d \text{Tr}(G_l^*G_l) a_l^2,\end{split}\end{equation*}
so that (\ref{Ho'}) holds for $\ll:= (1-\theta) \inf_{1\le l\le d} \|G_l\|_{HS}^2>0.$

A simple example such that  {\bf (A1)} holds is that $\si=I_{d\times d}, d=m-1$ and  $A_l$ given in (\ref{A}) with $\aa_l\ne\bb_l, 1\le l\le d.$
In this case we have $G_l^*G_k=0$ for $l\ne k$, so that {\bf (A1)} holds for $\theta =0.$

 The main tool in the study is the integration by parts formula of the Malliavin gradient. For fixed $T>0$, let $(D,\D(D))$ be the Malliavin gradient operator for the Brownian motion $\{B(t)\}_{t\in [0,T]}$, and let $(D^*,\D(D^*))$ be the adjoint operator. For any $F\in \D(D)$, the Malliavin gradient $DF$ is an element in $L^2(\OO\to\H;\P),$ where
$$\H:=\{\bb\in C([0,T];\R^d):\ \int_0^T|h'(t)|^2\d t<\infty\}$$ is the Carmeron-Martin space.
 We have
\beq\label{INT} \E \big[D_h \{f(X(T), Y(T))\}\big]= \E\big[f(X(T),Y(T)) D^*h\big],\ \ f\in C_b^1(\R^{m+d}), h\in\D(D^*).\end{equation}

To establish Bismut (resp. Driver) type formulae using (\ref{INT}), we need to construct element $h\in \D(D^*)$ such that the right-hand side of (\ref{INT}) reduces to $\nn_{(u,v)} P_T f$(resp.  $P_T\nn_{(u,v)} f$) for   given $(u,v)\in\R^{m+d}$. To this end,  let
 $Q_T=(q_{kl}(T))_{1\le l,k\le d}$ be a $\R^d\otimes\R^d$-valued random variable,
 where
 $$q_{lk}(T)=\int_0^T\bigg\<G_l^*G_k\Big(B(t)-\ff 1 T\int_0^TB(s)\d s\Big), B(t)-\ff 1 T\int_0^TB(s)\d s\bigg\>\d t,\ \ 1\le l,k\le d.$$
 Moreover,
 let $\aa_{T,u,v},\tt\aa_{T,u,v}\in\R^d$ with components
\beg{equation*}\beg{split} &(\aa_{T,u,v})_l = v_l -\<\si^{-1}u, A_l X(0)\> -\<A_l u, B(T)\> -\ff 1 T \int_0^T \<G_l^*\si^{-1}u, B(t)\>\d t,\\
&(\tt\aa_{T,u,v})_l = v_l -\<\si^{-1}u, A_l X(0)\> -\ff 1 T \int_0^T \<G_l^*\si^{-1}u, B(t)\>\d t,\ \ \ \  1\le l\le d.\end{split}\end{equation*}

\subsection{Main results} 

 \beg{thm}\label{T2.1} Assume {\bf (A1)} and let $T>0$ and $(u,v)\in \R^{m+d}$ be fixed. Let $h$ and $\tt h$ be such that $h(0)=\tt h(0)=0$ and
 \beg{equation*}\beg{split} & h'(t)= \ff{\si^{-1} u}T +\sum_{k=1}^d (Q_T^{-1}\aa_{T,u,v})_k G_k \bigg(B(t)-\ff 1 T\int_0^TB(s)\d s\bigg),\\
 &\tt h'(t)= \ff{\si^{-1} u}T +\sum_{k=1}^d (Q_T^{-1}\tt \aa_{T,u,v})_k G_k \bigg(B(t)-\ff 1 T\int_0^TB(s)\d s\bigg),\ \ t\in [0,T].\end{split}\end{equation*} Then:
 \beg{enumerate} \item[$(1)$] $h,\tt h\in \D(D^*)$ and for any $p>1$ there exists a constant $c_p>0$ independent of $(u,v)\in \R^{m+d}$ and $T>0$  such that
 $$\E |D^* h|^p + \E |D^*\tt h|^p\le \ff{c_p}{T^p}\big\{|v|^p +   |u|^p(|X(0)|^p +T^{\ff p 2 })\big\}.  $$
 \item[$(2)$] For any $f\in C_b^1(\R^{m+d}),$
 $P_T(\nn_{(u,v)} f)  =\E \big[f(X(T),Y(T)) D^*h\big].$
 \item[$(3)$] For any $f\in C_b^1(\R^{m+d})$,
 $\nn_{(u,v)}P_T f = \E \big[f(X(T),Y(T)) D^*\tt h\big].$\end{enumerate}Consequently,  for any $p>1$ there exists a constant $c_p>0$ such that
 $$\big(|P_T\nn_{(u,v)} f|+|\nn_{(u,v)}P_Tf|\big)(x,y) \le (P_T|f|^p)^{\ff 1 p} \ff{c_p}T \big\{|v|+|u|(|x|+\ss T\big)\big\}$$ holds for all $(u,v),(x,y)\in \R^{m+d}, T>0$ and $f\in C_b^1(\R^{m+d}).$ \end{thm}
 
 As a consequence of Theorem \ref{T2.1}, we have the following estimate (\ref{LB0}) of $\GG(P_tf)$, where
$$\GG(f):= \ff 1 2 \sum_{i=1}^m (X_if)^2$$ is the energy form associated to $L$. This estimate will imply the coupling property of the diffusion process as well as the Liouville property for the time-space harmonic functions. Recall that the $L$-diffusion process has the coupling property if for any initial points $z,z'\in\R^{m+d}$ one may construct two processes $Z_t,Z_t'$ generated by $L$   starting at $z,z'$ respectively, such that the coupling time $\tau:=\inf\{t\le 0: Z_t=Z_t'\}<\infty$. In this case $(Z_t,Z_t')$ is called a successful coupling of the process. Moreover, a bounded function $u$ on $[0,\infty)\times \R^{m+d}$ is called time-space harmonic associated to $P_t$, if $P_s u(t,\cdot)= u(t-s,\cdot)$ holds for any $t\ge s\ge 0.$ In particular, a bounded harmonic function is a time-space harmonic function.

Let  $\rr$ be the distance induced by $\GG$, i.e.
\beq\label{RH} \rr(z,z')= \sup\{|f(z)-f(z')|:\ f\in C^1(\R^{m+d}), \GG(f)\le 1\}.\end{equation}

\beg{cor}\label{C2.2} For any $p>1$ there exists a constant $c_p>0$ such that
\beq\label{LB0}\ss{\GG(P_tf)}\le \ff{c_p} {\ss t}  (P_t |f|^p)^{1/p},\ \ t>0, f\in \B_b(\R^{m+d}).\end{equation}
Consequently: \beg{enumerate} \item[$(1)$] Let  $P_t(z,\cdot)$ be the transition probability kernel of $P_t$, and let $\|\cdot\|_{var}$ be the totally variational  norm. There exists a constant $c>0$ such that
$$\|P_t(z,\cdot)-P_t(z',\cdot)\|_{var} \le \ff{c\rr(z,z')}{\ss t},\ \ t>0, z,z'\in\R^{m+d}.$$
\item[$(2)$] The $L$-diffusion process has the coupling property.
\item[$(3)$] Any time-space harmonic function associated to $P_t$ has to be constant. \end{enumerate}\end{cor}
 
 \subsection{Proofs}
To prove Theorem \ref{T2.1}(1), we need the following lemmas.

\beg{lem}\label{L2.2} Assume {\bf (A1)}. Then $Q_T$ is invertible and for any $p>1$ there exists a constant $C_p>0$ independent of $T>0$ such that
$$\E \|Q_T^{-1}\|^p \le \ff{C_p}{T^{2p}},\ \ T>0,$$ where $\|\cdot\|$ is the operator norm. \end{lem}

\beg{proof} Let $\bar Q_T=\text{diag}\, Q_T$; that is, $\bar Q_T= (q_{lk}(T) 1_{\{l\}}(k))_{1\le k,l\le d}.$ By {\bf (A1)},  $\bar Q_T$ is invertible and $Q_T\ge (1-\theta)\bar Q_T.$  Therefore, it suffices to show that
\beq\label{LE} \E q_{ll}(T)^{-p} \le \ff {C_p}{T^{2p}},\ \ T>0, 1\le l\le d\end{equation} holds for some constant $C_p>0$ independent of $T>0.$  Let
$e_l\in\R^d$ with $|e_l|=1$ such that $|G_l^*e_l|= \|G_l\|>0.$ Then
\beg{equation*}\beg{split} q_{ll}(T) &=\int_0^T \bigg|G_l\bigg(B(t)-\ff 1 T\int_0^T B(s)\d s\bigg)\bigg|^2\d t\\
&\ge \int_0^T \bigg\<G_l \bigg(B(t)-\ff 1 T \int_0^TB(s)\d s\bigg), e_l\bigg\>^2\d t\\
&= \|G_l\|^2 \int_0^T \bigg|b_l(t)- \ff 1 T \int_0^T b_l(s)\d s\bigg|^2 \d t,\end{split}\end{equation*} where
$$b_l(t) := \Big\<B(t), \ff{G_l^* e_l}{|G_l^* e_l|}\Big\>,\ \ t\ge 0$$ is an one-dimensional Brownian motion. Therefore,
$$\E\, q_{ll}(T)^{-p}\le \ff{1}{\|G_l\|^{2p}} \E\ff 1 {(\int_0^T |b_l(t)-\ff 1 T \int_0^T b_l(s)\d s |^2 \d t)^p}.$$
Combining this with
\beg{equation*}\beg{split} \int_0^T \bigg|b_l(t)-\ff 1 T \int_0^T b_l(s)\d s \bigg|^2 \d t&= \ff 1 {2T} \int_{[0,T]^2} |b_l(t)-b_l(s)|^2\d t\d s\\
&\ge \ff 1 {2T} \int_0^{\ff T 3}\d s\int_{\ff{2T}3}^T |b(t)-b(s)|^2\d t,\end{split}\end{equation*}
and using the Jensen inequality, we obtain
\beg{equation*}\beg{split} \E\, q_{ll}(T)^{-p}& \le \ff{6^{p}}{\|G_l\|^{2p}} \E\ff 1 {(\ff 3 T\int_0^{\ff T 3}\d s \int_{\ff {2 T}3}^T|b(t)-b(s)|^2\d t)^{p}}\\
&\le \ff{6^{p+1}}{\|G_l\|^{2p}T}\int_0^{\ff T 3} \E \Big(\ff 1 {(\int_{\ff {2T}3}^T |b(t)-b(s)|^2\d t)^{p}}\Big)  \d s.\end{split}\end{equation*}
According to \cite[Lemma 3.3]{W12a}, this implies (\ref{LE}) for some constant $C_p$ independent of $T>0$, and we thus finish the proof.
\end{proof}

\beg{lem}\label{L2.3} Assume {\bf (A1)}. Then $Q_T^{-1} \aa_{T,u,v}, Q_T^{-1}\tt \aa_{T,u,v}\in \D(D)^{\otimes d}$, and there exists a constant $c>0$ independent of $T>0$ such that for any adapted random variable   $\bb$ on the cameron-Martin space $\H$,
\beg{equation*}\beg{split} &|D_\bb Q_T^{-1}\aa_{T,u,v}|+ |D_\bb Q_T^{-1}\tt \aa_{T,u,v}|\\ &\le c T \|Q_T^{-1}\|^2\|\bb\|_\infty\|B\|_\infty\big\{|v|+|u|(|X(0)|+\|B\|_\infty)\big\}+ c\|Q_T^{-1}\|\cdot |u|\cdot\|\bb\|_\infty,\end{split}\end{equation*}where $\|\cdot\|_\infty$ is the uniform norm on $C([0,T];\R^d)$.
\end{lem}
\beg{proof} We only prove the desired upper bound for $\|D_\bb Q_T^{-1}\aa_{T,u,v}\|$, since that for the other term is completely similar. It is easy to see that
\beq\label{AA}|\aa_{T,u,v}|\le c_1 \big\{|v|+|u|(|X(0)|+\|B\|_\infty)\big\}\end{equation}  holds for some constant $c_1>0.$ Moreover, since $D_\bb B(t)= \bb(t), t\in [0,T]$, it follows from  the definitions of $q_{kl}(T)$ and $\aa_{T,u,v}$ that each components of $Q_T$ and $\aa_{T,u,v}$ are in $\D(D)$ with
$$|D_\bb q_{kl}(T)|\le c_2 \|\bb\|_\infty T\|B\|_\infty,\ \
 |D_\bb \aa_{T,u,v}|\le c_2 |u|\cdot\|\bb\|_\infty$$ holding for some constant $c_2>0$ and all $1\le k,l\le d.$  Combining these with the fact that
 $$D_\bb  Q_T^{-1}\aa_{T,u,v}= -Q_T^{-1} \{D_\bb Q_T \} Q_T^{-1}\aa_{T,u,v} + Q_T^{-1} D_\bb \aa_{T,u,v},$$ we derive the desired upper bound estimate of
$\|D_\bb Q_T^{-1}\aa_{T,u,v}\|$. \end{proof}

\beg{proof}[Proof of Theorem \ref{T2.1}] (1) We only prove for $h$ as that for $\tt h$ is similar. Let  $\{\ee_i\}_{1\le i\le m}$ be the canonical ONB of $\R^m$. Then $a_i:= \<a, \ee_i\>$ is the $i$-th coordinate of $a\in\R^m.$ Let
$$h_0(t)= \ff t T \si^{-1} u,\ \ h_i(t)= t\ee_i, \ \ \bb_k(t)= \int_0^t G_k B(s)\d s,\ \ 1\le i\le m, 1\le k\le d.$$ We have
$$h(t)= h_0(t)+ \sum_{k=1}^d \big(Q_T^{-1}\aa_{T,u,v}\big)_k \bb_k(t)- \sum_{i=1}^m \bigg(\sum_{k=1}^d\ff{(Q_t^{-1}\aa_{T,u,v})_k}T\int_0^T (G_kB(t))_i\d t\bigg)h_i(t)$$  and
\beg{equation*}\beg{split} &D^* h_0 =\ff 1 T \int_0^T \<\si^{-1} u, \d B(t)\>=\ff 1 T\<\si^{-1}u, B(T)\>,\ \ D^* h_i= B_i(T), \\
&D^*\bb_k = \int_0^T \<G_k B(t), \d B(t)\>, \ \   D_{h_i} B(t)= h_i(t)= t\ee_i,\ \ 1\le i\le m, 1\le k\le d. \end{split}\end{equation*}  Combining these with Lemma \ref{L2.3} and the fundamental identity
$$D^* (F\bb)= F D^*\bb - D_\bb F$$ for   $F\in \D(D), \bb\in\D(D^*)$ such that $F D^*\bb - D_\bb F\in L^2(\P)$,
we conclude that $h\in \D(D^*)$ and
\beg{equation*}\beg{split} D^* h=& \ff 1 T \<\si^{-1}u, B(T)\> +\sum_{k=1}^d  (Q_T^{-1}\aa_{T,u,v})_k \int_0^T\<G_k B(t), \d B(t)\> \\
&-\sum_{k=1}^d D_{\bb_k}(Q_T^{-1}\aa_{T,u,v})_k -\sum_{k=1}^d \sum_{i=1}^m \bigg(\ff{(Q_T^{-1}\aa_{T,u,v})_k}T \int_0^T (G_k B(t))_i\d t \bigg)B_i(T) \\
&+ \sum_{i=1}^m\sum_{k=1}^d  \ff{D_{h_i} (Q_T^{-1}\aa_{T,u,v})_k}T \int_0^T (G_k B(t))_i\d t
 + \sum_{i=1}^m\sum_{k=1}^d \ff   T 2 (Q_T^{-1}\aa_{T,u,v})_k (G_k)_{ii}.\end{split}\end{equation*}
Therefore, it is easy to see from  Lemma \ref{L2.3}, (\ref{AA}) and $(G_k)_{ii}=0$ that
\beg{equation*}\beg{split} |D^*h|\le & \ff{C|u|\cdot \|B\|_\infty}{T}   + C \|Q_T^{-1}\|\big\{|v|+|u|(|X(0)|+\|B\|_\infty\big\}\sum_{k=1}^d \Big|\int_0^T\<G_k B(t),\d B(t)\>\Big|\\
&\quad + CT^2\|Q_T^{-1}\|^2 \cdot\|B\|_\infty^2 \big\{|v|+|u|(|X(0)|+\|B\|_\infty)\big\}+ CT|u|\cdot\|Q_T^{-1}\|\cdot\|B\|_\infty  \\
&\quad +C \|Q_T^{-1}\|\cdot\|B\|_\infty^2 \big\{|v|+|u|(|X(0)|+\|B\|_\infty)\big\}\\
&\le  C|u|\cdot \|B\|_\infty\Big(\ff 1 T+ T\|Q_T^{-1}\|\Big) + C\|Q_T^{-1}\|\big\{|v|+|u|(|X(0)|+\|B\|_\infty)\big\}\\
&\qquad\qquad \qquad \times \Big(\sum_{l=1}^d \Big|\int_0^T\<G_k B(t),\d B(t)\>\Big|+T^2\|Q_T^{-1}\|\cdot \|B\|_\infty^2 +\|B\|_\infty^2\Big)\end{split}\end{equation*}
  holds for some constant $C>0.$ Combining this with Lemma \ref{L2.2} and the fact that for any $p>1$
 $$\E \|B\|_\infty^{2p}+ \sum_{k=1}^d \E \Big|\int_0^T \<G_k B(t), \d B(t)\>\Big|^p \le c(p) T^{p},\ \ T>0$$ holds for some constant $c(p)>0$, we obtain the desired upper bound of $\E |D^*h|^p.$

(2) For $\bb(s)= \sum_{i=1}^n \xi_i \bb_i(s)$, where $\xi_i$ are real-valued random variables and $\bb_i(s)$ are square-integrable adapted processes on $\R^d$, define
$$\int_0^t \<\bb(s),\d B(s)\>= \sum_{i=1}^n \xi_i \int_0^t \<\bb_i(s), \d B(s)\>.$$ Then it is easy to see from   (\ref{E1})     that
$$\beg{cases} D_h X(t)= \si h(t), \ D_hX(0)=0,\\
D_h Y_l(t)= \int_0^t \<A_l \si h(s), \d B(s)\> +\int_0^t \<A_l X(s), h'(s)\> \d s,\ D_h Y_l(0)=0, \ 1\le l\le d.\end{cases}$$ In particular,
\beq\label{P1} D_h X(T)= \si h(T)=u.\end{equation} Noting that   $X(s)= X(0)+ \si B(s)$ and $h(T)= \si^{-1} u$, we obtain
\beg{equation*}\beg{split} D_h Y_l(T) &= \<A_l \si h(T), B(T)\> -\int_0^T \<A_l\si h'(t), B(t)\>\d t +\int_0^T \<\si^* A_l^* h'(t), \si^{-1} (X(0) +\si B(t))\>\d t\\
&= \<A_l u, B(T)\> +\<A_l^* \si^{-1} u, X(0)\> +\int_0^T \<G_l^* h'(t), B(t)\>\d t.\end{split}\end{equation*} Combining this with the definition of $h'(t)$ and letting
$$\hat B(t)= B(t)-\ff 1 T \int_0^T B(s)\d s,\ \ t\in [0,T],$$ we arrive at
\beq\label{P2} \beg{split}D_h Y_l(T) &= \<A_l u, B(T)\> +\<A_l^* \si^{-1} u, X(0)\> \\
&\qquad +\int_0^T \Big\< \sum_{k=1}^d (Q_T^{-1}\aa_{T,u,v})_k   G_l^*G_k \hat B(t)  +\ff{G_l^*\si^{-1}u}T, B(t)\Big\>\d t\\
&= \<A_l u, B(T)\> +\<A_l^* \si^{-1} u, X(0)\>+\ff 1T \int_0^T \<G_l^*\si^{-1}u, B(t)\>\d t \\
&\qquad +\sum_{k=1}^d (Q_T^{-1}\aa_{T,u,v})_k \int_0^T \big\<G_l^*G_k \hat B(t), \hat B(t) \big\> \d t\\
&= \<A_l u, B(T)\> +\<A_l^* \si^{-1} u, X(0)\> +\ff 1 T \int_0^T \<G_l^*\si^{-1} u, B(t)\>\d t +(\aa_{T,u,v})_l=v_l.
\end{split}\end{equation}By (\ref{P1}) and (\ref{P2}) we obtain
$$D_h(X(T), Y(T))= (u,v).$$ Therefore, it follows from (\ref{INT}) that
\beg{equation*}\beg{split} P_T (\nn_{(u,v)} f) &= \E\big\<\nn f(X(T), Y(T)), (u,v)\big\> = \E\big\<\nn f(X(T), Y(T)), D_h(X(T), Y(T))\big\>\\
&= \E D_h \big\{f(X(T),Y(T))\big\} = \E\big\{f(X(T),Y(T)) D^* h\big\}.\end{split}\end{equation*}

(3)  Similarly to (2),  we have
$$\beg{cases} D_{\tt h} X(t)= \si \tt h(t), \\
D_{\tt h} Y(t)= \int_0^t \<A_l \si \tt h(s), \d B(s)\> +\int_0^t \<A_l X(s), \tt h'(s)\> \d s.\end{cases}$$ In particular,
\beq\label{P1'} D_{\tt h} X(T)= \si \tt h(T)=u.\end{equation} Noting that   $X(s)= X(0)+ \si B(s)$ and $\tt h(T)= \si^{-1} u$, as in (2) we obtain
$$ D_{\tt h} Y(T) = \<A_l u, B(T)\> +\<A_l^* \si^{-1} u, X(0)\> +\int_0^T \<G_l^* \tt h'(t), B(t)\>\d t.$$ Combining this with the definition of $\tt h'(t)$   we arrive at
\beq\label{P2'}\beg{split} D_{\tt h} Y_l(T)
&= \<A_l u, B(T)\> +\<A_l^* \si^{-1} u, X(0)\> +\ff 1 T \int_0^T \<G_l^*\si^{-1} u, B(t)\>\d t +(\tt \aa_{T,u,v})_l\\
&=v_l +  \<A_lu,   B(T)\>.\end{split}\end{equation}
  Moreover, it is easy to see that
$$\beg{cases} \d \nn_{(u,v)} X(t)= 0,\ \ \nn_{(u,v)}X(0)=u,\\
\d \nn_{(u,v)} Y_l(t)= \<A_l \nn_{(u,v)}X(t), \d B(t)\>,\ \ \nn_{(u,v)}Y_l(0)= v_l, \ \ 1\le l\le d.\end{cases}$$
Then
\beq\label{PP2}\nn_{(u,v)}X(T)=u,\ \ \nn_{(u,v)} Y_l(T)= v_l +\int_0^T \<A_l u, \d B(t)\>= v_l +\<A_l u, B(T)\>, \ \ 1\le l\le d.\end{equation} Combining this with
 (\ref{P1'}) and (\ref{P2'}) we obtain
$$D_{\tt h}(X(T), Y(T))= (\nn_{(u,v)}X(T), \nn_{(u,v)}Y(T)).$$ Therefore, it follows from (\ref{INT}) that
\beg{equation*}\beg{split} \nn_{(u,v)}P_T  f  &= \E\big\<\nn f(X(T), Y(T)), \nn_{(u,v)}(X(T),Y(T))\big\> = \E\big\<\nn f(X(T), Y(T)), D_{\tt h}(X(T), Y(T))\big\>\\
&= \E D_{\tt h} \big\{f(X(T),Y(T))\big\} = \E\big\{f(X(T),Y(T)) D^* \tt h\big\}.\end{split}\end{equation*}
\end{proof}

\beg{proof}[Proof of Corollary \ref{C2.2}]   By an approximation argument, it suffices to prove $(\ref{LB0})$ for $f\in C_0^\infty(\R^{m+d})$. Indeed, for any $z\in \R^{m+d}$, let $\ee\in\R^{m+d}$ be a unit vector such that $\ss{\GG(P_tf)}(z)= \nn_{\ee} P_tf(z).$ Since $P_t f\in C_b(\R^{m+d})$ for $f\in \B_b(\R^{m+d})$ and $t>0$,   (\ref{LB0}) holds at point $z$ provided
\beq\label{LB} \ff{P_t f(z+\vv v)- P_t f(v)}\vv \le \ff{c_p} {\ss t}   \ff 1 \vv \int_0^\vv (P_t |f|^p)^{1/p}(z+sv)\d s,\ \ \vv\in (0,1).\end{equation}    Noting that   (\ref{LB0})  with $f\in C_0^\infty(\R^{m+d})$ also implies  (\ref{LB})   for $f\in C_0^\infty(\R^{m+d})$, and that $C_0^\infty(\R^{m+d})$ is dense in $L^p\big(P_t(z,\cdot)+P_t(z+\vv \ee,\cdot)+\int_0^\vv P_t(z+s\ee,\cdot)\d s\big)$, we conclude that  (\ref{LB0})  for $f\in C_0^\infty(\R^{m+d})$ implies (\ref{LB}) for all $f\in \B_b(\R^{m+d})$, and hence also implies (\ref{LB0}) for all $f\in \B_b(\R^{m+d})$.

Next, by the left-invariant property of $X_i$,   it suffices to prove the desired estimate at   point $(0,0)\in \R^{m+d}$. To see this, for any $z\in\R^{m+d}$, let
$$\ell_z(z')= z\bullet z',\ \ z'\in\R^{m+d}.$$ Since $X_i$ are left-invariant, we have
$$\GG(P_t f)(z)= \GG((P_tf)\circ\ell_z)(0,0)= \GG(P_tf\circ\ell_z)(0,0),$$ so that the desired estimate at point $(0,0)$ implies
$$ \GG(P_t f)(z)\le \ff{c_p} {\ss t }  (P_t |f\circ\ell_z|^p)^{1/p}(0,0) = \ff{c_p} {\ss t }  (P_t |f|^p)^{1/p}(z).$$

Now, we intend to prove (\ref{LB0})  for $f\in C_0^\infty(\R^{m+d})$ at point $(0,0)$. In this case, there exists an unit element $u\in\R^m$ such that $\ss{\GG(P_tf)}(0,0)= \nn_{(u,0)}P_tf(0,0)$. Then, by Theorem \ref{T2.1} (1) and (3) with $X(0)=0$ and using the H\"older inequality, we derive the desired upper bound  for $\ss{\GG(P_tf)}(0,0).$

Finally, noting that $$\|P_t(z,\cdot)-P_t(z',\cdot)\|_{var}=2\sup_{\|f\|_\infty\le 1} |P_tf(x)-P_t f(y)|\le 2\rr(x,y) \sup_{\|f\|_\infty\le 1} \ss{\|\GG(P_tf)\|_\infty},$$ then (\ref{LB0}) implies (1). According to \cite{CG} (see also \cite{CW}), (2) and (3) follow from (1).\end{proof}

\section{An explicit   inverse Poincar\'e inequality} Note that by Corollary \ref{C2.2},  {\bf (A1)} implies 
\beq\label{GD}\GG(P_tf)\le \ff C t \big(P_t f^2- (P_Tf)^2\big),\ \ f\in\B_b(\R^{d+m}), t>0\end{equation} for some constant $C>0$. In fact, this estimate follows also from \eqref{Ho} according to \cite[Proposition 4.7]{BGM}.
In this  section we aim  to prove this inequality with   an explicit $L^2$-estimate  on $\GG (P_t f)$ as  in \cite[Section 3]{BBBC}, where the heat semigroup on the Heisenberg group is concerned. To this end, we need the following assumption:\beg{enumerate} \item[{\bf (A2)}]
 For any $l,l',l''\in\{1,\cdots,d\}$, $ A^*_l= -A_l, \si A_l= A_l \si, A_l A_{l'}= A_{l'}A_l,$ and $A_l\si, A_l A_{l'}A_{l''}\si$ and $A_l\si^2\si^*$ are skew-symmetric. \end{enumerate}
A simple example for this assumption to hold is that $\si=I_{m\times m}$ and $\{A_l\}$ are commutative skew-symmetric $m\times m$-matrices.

\beg{thm}\label{T4.1} Assume $(\ref{Ho})$ and  {\bf (A2)}. Then for any $f\in \B_b(\R^{m+d})$ and $t>0,$
$$\GG(P_t f)\le \ff{m+2d}{2t} \{P_t f^2-(P_t f)^2\}.$$ \end{thm}

This estimate is equivalent to (\ref{LB0}) for $p=2$ with explicit constant $c_p=\big(\ff{m+2d}2\big)^{\ff 1 2}.$
To prove this result,  we introduce the dilation operator modified from \cite{BBBC},
$$\DDD:=\ff 1 2 \sum_{i=1}^mx_i \pp_{x_i} +\sum_{l=1}^d y_l\pp_{y_l}$$ and the dual vector fields
$$\hat X_i(x,y):= X_i(x,y) - 2 \sum_{l=1}^d (A_l x)_i \pp_{y_l},\ \ \ 1\le i\le m.$$
Simply denote
$$p_t^0(z)= p_t((0,0), z),\ \ z\in\R^{m+d},$$ where $p_t$ is the heat kernel of $P_t$ w.r.t. the Lebesgue measure $\mu$ on $\R^{m+d},$ which exists due to the H\"ormander condition.

\beg{lem}\label{L4.2} Assume $(\ref{Ho})$ and  {\bf (A2)}. Then $L\DDD -\DDD L =L$ and $[\hat X_i, X_j]=0,\ 1\le i,j\le m$. Consequently, $(tL+\DDD +\ff {m+2d}2) p_t^0  =0$ and
$\hat X_i P_t= P_t \hat X_i,\ 1\le i\le m.$ \end{lem}

\beg{proof} It is easy to see that $[X_i,\DDD] =\ff 1 2 X_i, 1\le i\le m.$ Then
\beq\label{DB} L\DDD-\DDD L= \ff 1 2 \sum_{i=1}^m (X_i^2\DDD -\DDD X_i^2)=\ff 1 2 \sum_{i=1}^m (X_i[X_i,\DDD]+[X_i,\DDD] X_i) =\ff 1 2\sum_{i=1}^m X_i^2 =L.\end{equation} Let $T_t$ be the semigroup generated by $\DDD$. Then $T_t\e^t $ is generated by $\DDD+1$  and due to (\ref{DB})
$ L \DDD= (\DDD+1)L$. So, $L T_s= T_s \e^s L,$ which implies that
\beq\label{DB1}  P_t T_s= T_s P_{\e^{s} t},\ \ t,s\ge 0.\end{equation}  Differentiating both sides w.r.t. $s$ at $s=0$, we obtain
\beq\label{DB2} P_t\DDD = \DDD P_t +t P_t L,\ \ t\ge 0.\end{equation}
Since $\DDD(0,0)=0$, it follows that
$$P_t (tL-\DDD)f(0,0)=0,\ \ f\in C_0^\infty(\R^{m+d}).$$ Combining this with
\beg{equation*} \beg{split}P_t(tL-\DDD)f(0,0)&= \int_{\R^{m+d}} p_t^0(z) (tL-\DDD) f(z)\d z\\
& = \int_{\R^{m+d}} f(z) \big\{(tL+\DDD)p_t^0(z) + (\text{div}\DDD)p_t^0(z)\big\}\d z\\
&=\int_{\R^{m+d}} f(z) \Big(tL+\DDD+\ff {m+2d}2\Big)p_t^0 (z) \d z,\end{split}\end{equation*} we conclude that $(tL+\DDD+\ff {m+2}2)p_t^0 =0.$

Next, for any $1\le i,j\le m$, we have \beg{equation*}\beg{split} &[\hat X_i, X_j]= [X_i,X_j] + 2 \sum_{k=1}^m\sum_{l=1}^d(A_l)_{ik}\si_{kj}\pp_{y_l} \\
&= \sum_{l=1}^d \big\{(A_l\si)_{ji}-(A_l\si)_{ij}\big\}\pp_{y_l} + 2\sum_{l=1}^d (A_l\si)_{ij}\pp_{y_l}=0\end{split}\end{equation*} since $A_l\si$ is skew-symmetric. This implies $\hat X_i P_t= P_t \hat X_i$ for any $1\le i\le m$.
\end{proof}

\beg{lem}\label{L4.3} Assume  {\bf (A2)} and let $\hat \GG(f)= \ff 1 2 \sum_{i=1}^m (\hat X_i f)^2$. Then
$\hat\GG(p_t^0) =\GG(p_t^0).$ \end{lem}

\beg{proof} It is easy to see that at point $(x,y)\in\R^{m+d},$
\beg{equation*}\beg{split} \hat\GG(f) &= \GG(f) - 2\sum_{i=1}^m (X_i f) \Big(\sum_{l=1}^d (A_lx)_i\pp_{y_l}  f\Big) + 2 \sum_{i=1}^m\Big(\sum_{l =1}^d (A_lx)_i  \pp_{y_l}  f\Big)^2\\
&=\GG(f) - 2\sum_{l=1}^d \sum_{k=1}^m (\si A_l x)_k \pp_{x_k} f)(\pp_{y_l}f)=\GG(f) - 2\sum_{l=1}^d (\pp_{y_l}f) \Theta_lf,\end{split}\end{equation*} where
$$\Theta_l:= \sum_{k=1}^m (\si A_l x)_k\pp_{x_k},\ \ 1\le l\le d.$$ So, it remains to prove $\Theta_l p_t^0 =0$ for $1\le l\le d.$  We prove it by two steps.

(1) $\Theta_l L= L\Theta_l.$ Since $A_l\si=\si A_l$, it is easy to see that
\beg{equation*}\beg{split} [\Theta_l, X_i] &= \sum_{l'=1}^d \sum_{k=1}^m (\si A_l x)_k (A_{l'})_{ik}\pp_{y_{l'}} -\sum_{k,j=1}^m\si_{ji}(\si A_l)_{kj} \pp_{x_k}\\
&= \sum_{l'=1}^d (A_{l'}\si A_lx)_i\pp_{y_{l'}} -\sum_{k=1}^m (\si A_l \si)_{ki} \pp_{x_k} \\
&= \sum_{l'=1}^d (A_{l'} A_l\si x)_i \pp_{y_{l'}} -\sum_{i=1}^m (A_l \si^2)_{ki} \pp_{x_k}.\end{split}\end{equation*} Then
\beg{equation*}\beg{split} &\sum_{i=1}^m (\Theta_l X_i^2 -X_i^2 \Theta_l) = \sum_{i=1}^m \big\{[\Theta_l, X_i] X_i +X_i [\Theta_l, X_i]\big\}\\
&= 2 \sum_{l'=1}^d \sum_{i,k=1}^m (A_{l'}A_l \si x)_i \si_{ki} \pp_{x_k} \pp_{y_{l'}} + 2 \sum_{l',l''=1}^d \sum_{i=1}^m (A_{l'}A_l\si x)_i (A_{l''}x)_i \pp_{y_{l'}}\pp_{y_{l''}}\\
&\quad -2 \sum_{i,k,j=1}^m (A_l \si^2)_{ki} \si_{ji} \pp_{x_k}\pp_{x_j} - 2\sum_{i,k=1}^m\sum_{l'=1}^d (A_l \si^2)_{ki} (A_{l'} x)_i \pp_{x_k}\pp_{y_{l'}}\\
&\quad +\sum_{l'=1}^d \sum_{i,k=1}^m \si_{ki} (A_{l'}A_l\si)_{ik}\pp_{y_{l'}} -\sum_{l'=1}^d \sum_{i,k=1}^m(A_l\si^2)_{ki} (A_{l'})_{ik}\pp_{y_{l'}} \\
&= 2\sum_{l'=1}^d \sum_{k=1}^m \big\{(\si A_{l'}A_l\si x)_k -(A_l \si^2 A_{l'}x)_k\big\} \pp_{x_k}\pp_{y_{l'}}
+2 \sum_{l',l''=1}^d \big\<A_{l'}A_l\si x, A_{l''} x\big\> \pp_{y_{l'}} \pp_{y_{l''}}\\
&\quad -2 \sum_{i,j=1}^m (A_l \si^2 \si^*)_{ij}\pp_{x_i}\pp_{x_j} + \sum_{l'=1}^d \text{Tr} (\si A_{l'}A_l \si -A_{l'}A_l \si^2)\pp_{y_{l'}}.\end{split}\end{equation*} Due to {\bf (A2)}, this implies that $\Theta_l L= L\Theta_l.$

(2) By (1), $\text{div}\Theta_l=0$ and $ \Theta_l(0,0)=0$, for any $f\in C_0^\infty(\R^{m+d})$ we have
\beg{equation*}\beg{split} &0= (\Theta_l P_t f)(0,0) = (P_t\Theta_lf)(0,0) \\
&=\int_{\R^{m+d}} p_t^0(z) \Theta_l f(z)\d z = -\int_{\R^{m+d}} \{\Theta_l p_t^0(z)\}f(z)\d z.\end{split}\end{equation*} Therefore, $\Theta_l p_t^0 =0.$
\end{proof}

\beg{lem}\label{L4*} Assume $(\ref{Ho})$. Then there exists two constants $c_1,c_2>0$ such that
$$p_t(z,z')\le \ff{c_1\exp[-\ff{c_2\rr(z,z')^2}t]}{t^{(m+2 d)/2}},\ \ t>0, z,z'\in\R^{m+d}.$$\end{lem}

\beg{proof} We shall use the dimension-free Harnack inequality derived in \cite{BGM} using the generalized curvature condition. 
Let
$$\GG(f,g)= \ff 12\sum_{i=1}^m (X_if)(X_ig),\ \ \GG^Z(f,g)= \ff 12\sum_{l=1}^d (\pp_{y_l} f)(\pp_{y_l}g),\ \ f,g\in C^1(\R^{m+d})
$$ and denote $\GG(f)=\GG(f,f),\GG^Z(f)=\GG^Z(f,f)$. Define
$$\GG_2(f)= \ff 1 2 L \GG(f,f) -\GG(f, Lf),\ \ \GG_2^Z(f)= \ff 1 2 L\GG^Z(f,f) -\GG^Z(f, Lf),\ \ f\in C^3(\R^{m+d}).$$
By \cite[Proposition 4.4]{BGM}, \eqref{Ho} (equivalently, \eqref{Ho'}) implies the generalized curvature condition
$$\GG_2(f)+r \GG_2^Z(f)\ge c \GG^Z(f)-\ff{c'}r \GG(f),\ \ f\in C^2(\R^{m+d}), r>0$$ for some constants $c,c'>0,$ see also Lemma \ref{CD} below for a generalized curvature-dimension condition. According to \cite{BB1} (see also \cite[Proposition 4.7]{BGM}), this implies the following Harnack inequality of type \cite{W97} for some constant $C>0$:
$$(P_t f(z))^p\le(P_t f^p)(z')\exp\bigg[\ff{Cp}{(p-1)t} \rr(z,z')^2\bigg],\ \ t>0, z,z'\in\R^{m+d}, f\in\B_b^+(\R^{m+d}).$$  According to \cite{GW}, this Harnack inequality implies 
$$p_t(z,z')\le \ff{c_1\exp[-\ff{c_2\rr(z,z')^2}t]}{\ss{\mu(B(z, \ss t))\mu(B(z', \ss t))}},\ \ t>0, z,z'\in\R^{m+d} $$ for some constant $c_1,c_2>0$, where $\mu$ is the Lebesgue measure and  $B(z,r)=\{\rr(z,\cdot)\le r\}$ for $z\in\R^{m+d}$ and $r\ge 0.$  Since both $\rr$ and $\mu$ are left-invariant under the group action, this is equivalent to
$$ p_t(z,z')\le \ff{c_1\exp[-\ff{c_2\rr(z,z')^2}t]}{ \mu(B((0,0), \ss t))},\ \ t>0, z,z'\in\R^{m+d}.  $$
So, it remains to show that
\beq\label{MU} \mu(B((0,0), r))\ge c r^{m+2d},\ \ r>0\end{equation} holds for some constant $c>0.$ To see this, let us observe that for any $f\in C^1(\R^{m+d})$ and $f_r(x,y):=f(r x, r^2 y)$ one has
$$\GG(f_r)(x,y)= r^2  \GG(f)(rx, r^2y).$$ Combining this with with (\ref{RH}), we obtain $\rr((r x, r^2 y),(0,0))= r \rr((x,y), (0,0)).$ So,
\beg{equation*}\beg{split} B((0,0),r) &= \big\{(x,y)\in\R^{m+d}: \rr((x/r, y/r^2), (0,0))\le 1\big\}\\
&\supset \Big\{rx: \rr((x,0),(0,0)\le \ff 1 2\Big\}\times \Big\{r^2 y: \rr((0,y),(0,0))\le \ff 1 2\Big\}.\end{split}\end{equation*} Therefore, (\ref{MU}) holds for
$$c:= \mu\Big(\Big\{x: \rr((x,0),(0,0))\le \ff 1 2\Big\}\times \Big\{y: \rr((0,y),(0,0))\le \ff 1 2\Big\}\Big)\ge \mu\big(B((0,0), 1/2)\big)>0.$$
\end{proof}

\beg{lem}\label{L4.4} Assume $(\ref{Ho})$ and {\bf (A2)}. Then $\int_{\R^{m+d}} \GG(\log p_t^0, p_t^0)(z) \d z= \ff{m+2d}{2t},\ t>0.$\end{lem}

\beg{proof} We shall prove the lemma by using an approximation argument.  Let $h\in C_0^\infty([0,\infty))$ such that $0\le h\le 1, h|_{[0,1]}=1$ and $h|_{[2,\infty)}=0.$ Let $f_n(z)= h(z/n), n\ge 1.$ Then there is a constant $C_1>0$ such that
\beq\label{S1} |Lf_n|(z) + \GG(f_n)(z) +|\DDD f_n|(z) \le C_1 1_{[n,2n]}(|z|),\ \ z\in\R^{m+d}.\end{equation} Moreover, there is a constant $C_2>0$ such that $\GG(C_2\log (1+|\cdot|))\le 1.$ So, according to (\ref{RH})
$$\rr(0,z)\ge C_2\log (1+|z|),\ \ z\in \R^{m+d}.$$ Combining this with Lemma \ref{L4*} we obtain
\beq\label{S2} p_t^0(z) \le c_1(t) \exp\big[-c_2(t)\{\log (1+|z|)\}^2\big],\ \ z\in\R^{m+d}\end{equation} for some constants $c_1(t), c_2(t)>0.$
Since by Lemma \ref{L4.2} $(tL+\DDD+\text{div}\DDD)p_t^0=0$, for any $n\ge 1$ we have
\beg{equation*}\beg{split}  &\int_{\R^{m+d}} \big\{f_n\GG(\log p_t^0, p_t^0)\big\}(z)\d z = -t \int_{\R^{m+d}} \big\{f_n (\log p_t^0)L p_t^0 + (\log p_t^0)\GG(f_n, p_t^0)\big\}(z)\d z\\
&= \int_{\R^{m+d}} \big\{(f_n\log p_t^0)(\DDD+\text{div}\DDD)p_t^0\big\}(z)\d z- t\int_{\R^{m+d} } (\log p_t^0)(z) \GG(f_n, p_t^0) (z)\d z\\
&= -\int_{\R^{m+d}} \DDD(f_n p_t^0)(z)\d z +\int_{\R^{m+d}}\big\{(\DDD f_n) p_t^0 -(p_t^0\log p_t^0)\DDD f_n - t (\log p_t^0)\GG(f_n, p_t^0)\big\}(z)\d z\\
&= \ff{m+2d}2 \int_{\R^{m+d}} (f_np_t^0)(z)\d z +\int_{\R^{m+d}}\big\{(\DDD f_n) p_t^0 -(p_t^0\log p_t^0)\DDD f_n - t (\log p_t^0)\GG(f_n, p_t^0)\big\}(z)\d z\\
&\le \ff{m+2d} 2 + C(t) \int_{\{n\le |z|\le 2n\}}\Big\{p_t^0 +|p_t^0\log p_t^0|+|\log p_t^0|\ss{\GG(p_t^0)}\Big\}(z)\d z\end{split}\end{equation*}
for some constant $C(t)>0$ according to (\ref{S1}). Therefore, it suffices to verify
\beq\label{S3} \int_{\R^{m+d}}\Big\{p_t^0 +|p_t^0\log p_t^0|+|\log p_t^0|\ss{\GG(p_t^0)}\Big\}(z)\d z<\infty\end{equation} so that the desired estimate follows by letting $n\to\infty$. Noting that $p_t^0= P_{\ff t 2}p_{\ff t 2}^0$, \eqref{GD} and Lemma \ref{L4*} yield
$$\ss{\GG(p_t^0)}\le C_1(t) \ss{P_{\ff t 2} (p_{\ff t 2}^0)^2}\le C_2(t) \ss{p_t^0}$$ for some constants $C_1(t), C_2(t)>0.$ Therefore, (\ref{S3}) follows from (\ref{S2}) since
$$p_t^0+|p_t^0\log p_t^0|+ C_2(t)|\log p_t^0|\ss{p_t^0}\le C_3(t) \big\{(p_t^0)^2 + (p_t^0)^{\ff 1 4}\big\}$$ holds for some constant $C_3(t)>0.$  \end{proof}

\beg{proof}[Proof of Theorem \ref{T4.1}] As explained in \cite[Proof of Theorem 3.1]{BBBC}, we have
$$\GG(P_t f) \le \big\{P_t f^2-(P_tf)^2\big\}\int_{\R^{m+d}}\hat\GG(\log p_t^0, p_t^0)(z)\d z.$$ Then the proof is finished by combing this with Lemmas \ref{L4.3} and \ref{L4.4}.\end{proof}

\section{The Poincar\'e inequality}
In this section we prove the     estimate (\ref{GG}) below by following the argument in \cite[Section 4]{BBBC}. This estimate  for the heat semigroup on the Heisenberg group was first derived  in \cite{DM}.

According to (\ref{Ho}), there exists $\{(i_l,j_l)\}_{1\le l\le d}$ with   $1\le i_l<j_l\le m$ such that the matrix
 $$\tt M:= (\tt M_{l, l'})_{1\le l,l'\le d}$$ is invertible, where $\tt M_{l,l'}:= M_{(i_l,j_l),l}= (G_{l'})_{j_li_j}.$ Recall that for any $x\in\R^m$ and $1\le i\le m,$  $(A_\cdot x)_i:=((A_lx)_i)_{1\le l\le d} \in \R^d.$   Similarly, we let $(A_\cdot\si)_{ij}=((A_l\si)_{ij})_{1\le l\le d}\in\R^d.$

\beg{thm}\label{T5.1} Assume {\bf (A2)}. Then
\beq\label{GG}\GG(P_t f)\le C P_t \GG(f),\ \ t\ge 0, f\in C_b^1(\R^{m+d})\end{equation} holds for
$$C:= 2+ 16\sum_{i,j=1}^m|(\tt M^*)^{-1}(A_\cdot\si)_{ij}|^2 +32 P_1 \Big\{\sup_{1\le l\le d} \big\{(\tt M^*)^{-1}(A_\cdot x)_i\big\}_l^2\, \GG(\log p_1^0)\Big\}(0,0)<\infty,$$
where $p_1^0(z):= p_1((0,0),z)$ and $(A_\cdot x)_i(z):= (A_\cdot x')_i$ for $z=(x',y')\in\R^{m+d}.$   Consequently, the Poincar\'e inequality
$$P_t f^2 -(P_t f)^2\le 2C t P_t\GG(f),\ \ f\in C_b^1(\R^{m+d}) $$ holds for all $t>0.$
\end{thm}

To prove this result, we need the following lemma on curvature-dimension condition. When $m=\infty$ it reduces to the generalized curvature condition derived in \cite[Proposition 4.4]{BGM}. 

\beg{lem}\label{CD}   Assume the H\"ormander condition $(\ref{Ho'})$. For any $f\in C^3(\R^{m+d})$ and $r>0$,
$$\GG_2(f)+r \GG_2^Z(f) \ge \ff{(Lf)^2} m +\ff {c_2(G) \GG^Z(f)} 4 - \ff{c_1(G)}{r}\GG(f).$$
\end{lem}
\beg{proof} Recall that $L= \ff 1 2 \sum_{i=1}^m X_i^2$ and
$$[X_i,X_j]= -\sum_{l=1}^d (G_l)_{ij}\pp_{y_l},\ \ [X_i, \pp_{y_l}]=0,\ \ 1\le i,j\le m, 1\le l\le d.$$
Then
\beg{equation}\label{GG1}\beg{split} \GG_2(f) &= \ff 1 8 \sum_{i,j=1}^m X_i^2(X_j f)^2 -\ff 1 4 \sum_{j=1}^m \Big(X_j\sum_{i=1}^m X_i^2f\Big)(X_jf)\\
&=\ff 1 4 \sum_{i,j=1}^m (X_jf) (X_i^2X_j f)+\ff 1 4 \sum_{i,j=1}^m (X_iX_jf)^2 -\ff 1 4 \sum_{i,j=1}^m (X_jX_i^2f)(X_jf)\\
&= \ff 1 4 \sum_{i,j=1}^m (X_jf)([X_i,X_j]X_i+X_i[X_i,X_j])f +\ff 1 4 \sum_{i,j=1}^m (X_iX_jf)^2\\
&=\ff 1 4 \sum_{i,j=1}^m (X_iX_jf)^2-\ff 1 2 \sum_{i,j=1}^m (X_j f) \Big(\sum_{l=1}^d (G_l)_{ij} \pp_{y_l}X_i f\Big).\end{split}\end{equation}
Moreover,
$$ \GG_2^Z(f) = \ff 1 8 \sum_{l=1}^d \sum_{i=1}^m X_i^2 (\pp_{y_l} f)^2-\ff 1 4 \sum_{i=1}^m\sum_{l=1}^d (\pp_{y_l}f)(X_i^2\pp_{y_l}f)= \ff 1 4  \sum_{i=1}^m\sum_{l=1}^d (\pp_{y_l}X_i f)^2.$$ Combining this with (\ref{GG1}) and the fact
$$\ff 1 4 \sum_{i,j=1}^m (X_iX_jf)^2 \ge \ff {(Lf)^2} m +\ff 1 4 \sum_{1\le i\ne j\le m} (X_iX_jf)^2,$$ we obtain
\beq\label{GG2} \beg{split}\GG_2(f)+r\GG_2^Z(f) &\ge \ff {(Lf)^2} m +\ff 1 4 \sum_{1\le i\ne j\le m} (X_iX_jf)^2-\ff 1 {4r} \sum_{i=1}^m\Big(\sum_{j=1}^m\sum_{l=1}^d (X_jf)(G_l)_{ij}\Big)^2\\
&\ge  \ff {(Lf)^2} m +\ff 1 4 \sum_{1\le i\ne j\le m} (X_iX_jf)^2-\ff {c_1(G)} {2r} \GG(f,f),\ \ r>0.\end{split}\end{equation}
Finally, as observed in \cite{BG} we have
\beg{equation*}\beg{split} \sum_{1\le i\ne j\le m} (X_iX_jf)^2 &=\sum_{1\le i<j\le m}\big\{(X_iX_jf)^2+(X_jX_if)^2\big\}\\&=\ff 1 2\sum_{1\le i<j\le m} \big\{(X_iX_jf+X_jX_if)^2+(X_iX_jf-X_jX_if)^2\big\}\\
&\ge \ff 1 2 \sum_{1\le i<j\le m} ([X_i,X_j]f)^2 =\ff 1 2 \sum_{1\le i<j\le m} \Big(\sum_{l=1}^d (G_l)_{ij}\pp_{y_l} f\Big)^2\ge c_2(G) \GG^Z(f).\end{split}\end{equation*}
Combining this with (\ref{GG2}) we complete the proof.
\end{proof}As it is easy to see that the commutation condition
$$\GG(f, \GG^Z(f,f))=\GG^Z(f,\GG(f,f)),\ \ f\in C^\infty(\R^{m+d})$$ holds, the following assertions follow from the curvature-dimension condition presented in Lemma \ref{CD}, where  (4) and (5) are known as Li-Yau type gradient estimate and parabolic Harnack inequality (see \cite{LY}), and (2) is the dimension-free Harnack inequality initiated by the author in \cite{W97}, which implies the log-Harnack inequality (3) as observed in \cite{W10}. This type of Hanrack inequality was also established in \cite{DG} on a class of Lie groups. The entropy gradient inequality (1) implying the dimension-free Harnack inequality (2) was first observed in \cite{ATW06}.

\beg{cor}\label{C3.2} Assume the H\"ormander condition $(\ref{Ho'})$. For any $t>0$ and positive $f \in \B_b(\R^{m+d})$, the following assertions hold:
\beg{enumerate}\item[$(1)$] $\dfrac{t\GG(P_tf)}{P_t f}+ \dfrac{c_2(G)t^2 \GG^Z(P_tf)}{4P_tf}\le \dfrac{c_2(G)+8c_1(G)}{c_2(G)} \big\{P_t(f\log f)- (P_tf)\log P_t f\big\}.$
\item[$(2)$] $(P_t f)^p(z) \le (P_t f^p(z'))\exp\Big[\dfrac  {p(c_2(G)+ 8c_1(G))}{4(p-1)c_2(G)t}\rr(z,z')^2\Big],\ p>1, z,z'\in\R^{m+d}.$
\item[$(3)$] $P_t\log f(z)\le \log P_tf(z') + \dfrac{c_2(G)+ 8c_1(G)}{4c_2(G)t}\rr(z,z')^2,\ z,z'\in\R^{m+d}.$
\item[$(4)$] $\GG(\log P_t f) +\dfrac{c_2(G)t}{6} \GG^Z(\log P_tf) \le \dfrac{c_2(G)+6c_1(G)}{c_2(G)} \pp_t \log P_t f +\dfrac{m(c_2(G)+6c_1(G))^2}{2c_2(G)^2t}.$
\item[$(5)$] $P_t f(z)\le P_{t+s}f(z') \Big(\dfrac{t+s}t\Big)^{\ff{m(c_2(G)+6c_1(G))}{2c_2(G)}} \exp\Big[\dfrac{(c_2(G)+6c_1(G))\rr(x,y)^2}{4mc_2(G)s}\Big],\ \ s,t>0.$
\end{enumerate}
\end{cor}

\beg{proof}  According to Lemma \ref{CD},   \cite[Propositions 3.1, 3.4]{BB1} and \cite[Theorems 6.1, 7.1, 8.1]{BG},   it suffices to verify the following conditions:

(i) There exists a sequence $\{h_n\}\subset C_0^\infty(\R^{m+d})$ such that $h_n\uparrow 1$ and $\|\GG(h_n)\|_\infty+\|\GG^Z(h_n)\|_\infty\to 0$ as $n\uparrow\infty.$

(ii) $\GG(f,\GG^Z(f),f)= \GG^Z(f, \GG(f)),\ \ f\in C^\infty(\R^{m+d}).$

(iii) For any $f\in C_0^\infty(\R^{m+d})$ and $T>0$,
$$\sup_{t\in [0,T]} \big(\|\GG(P_t f)\|_\infty +\|\GG^Z(P_tf)\|_\infty\big)<\infty.$$

Let $f\in C_0^\infty([0,\infty))$ with $f'\le 0, f|_{[0,1]}=1$ and $f|_{[2,\infty)}=0$. Then (i) holds for $h_n(z):= f(|z|/n),\ n\ge 1, z\in\R^{m+d}.$ Next, (ii) follows from  $[X_i,\pp_{y_l}]=0,\ 1\le i\le m, 1\le l\le d.$  Finally, it is easy to see that Lemma \ref{CD} implies  assumption {\bf (A)} in \cite{WN} with $l=1, \GG^{(1)}=\GG^Z,  K_0(r)= -\ff{c_1(G)}r, K_1(r)= \ff{c_2(G)}4,$ and $W(x,y)= 1 +|x|^2 +|y|^2.$   Therefore,   (iii) is ensured by \cite[Lemma 2.1]{WN}, see also \cite[Lemma 5.2.2]{WB}.
\end{proof}

\beg{proof}[Proof of Theorem \ref{T5.1}] The desired Poincar\'e inequality follows immediately from (\ref{GG}) by noting that for $f\in C_0^\infty(\R^{m+d}),$
$$\ff{\d}{\d s} P_s (P_{t-s}f)^2 = 2 P_s \GG(P_{t-s}f) \le 2CP_t \GG(f),\ \ s\in [0,t].$$
Below we prove (\ref{GG}) and the finite of $C$ respectively, where the proof of (\ref {GG}) is modified from \cite{BBBC}.

(1) We first observe that   to prove (\ref{GG})  it suffices to confirm
\beq\label{GG'} \GG (P_1 f)(0,0) \le C P_1 \GG(f)(0,0),\ \ f\in C_0^\infty(\R^{m+d}).\end{equation}
Indeed, by the left-invariant property of $\GG$ and $P_t$, we only need to prove (\ref{GG}) at point $(0,0)$; and by a standard approximation argument as in the proof of Corollary \ref{C2.2}, we may assume that $f\in C_0^\infty(\R^{m+d})$. Finally, for any $t>0$, it follows from (\ref{DB1}) that
$$P_t f= P_t T_{-\log t}T_{\log t}= T_{-\log t}P_1 T_{\log t} f.$$ Noting that $T_sf(x,y) = f(\e^{\ff s 2}x, \e^s y)$, we have
$X_i T_s = \e^{\ff s 2} T_s X_i.$ Therefore, if (\ref{GG'}) holds, then at point $(0,0)$ we have
 \beg{equation*} \beg{split}  \GG(P_t f) &= \GG(T_{-\log t}P_1 T_{\log t} f) =\ff 1 t T_{-\log t} \GG(P_1T_{\log t}f)\\
  &\le \ff C t T_{-\log t} P_1\GG(T_{\log t} f) = C T_{-\log t} P_1 T_{\log t}\GG(f) = C P_t\GG(f).\end{split}\end{equation*}

(2) Note that
\beg{equation*}\beg{split} &\sum_{l=1}^d \big\{(\tt M^*)^{-1}(A_\cdot x)_i\big\}_l [X_{i_l},X_{j_l}] = -\sum_{l,l'=1}^d (\tt M^*)^{-1}_{ll'}(A_{l'}x)_i \tt M_{ll''}\pp_{y_{l''}}\\
&=-\sum_{l',l''=1}^d \big\{\tt M^*(\tt M^*)^{-1}\big\}_{l''l'} (A_{l'}x)_i \pp_{y_{l''}}=-\sum_{l'=1}^d (A_{l'}x)_i \pp_{y_{l'}}.\end{split}\end{equation*}  Then, by Lemma \ref{L4.2},  at point $(0,0)$ we have
\beq\label{PW}\beg{split}  X_i P_1 f &=\hat X_i P_1 f= P_1 \hat X_i f = P_1\Big\{X_i - 2 \sum_{l=1}^d (A_l x)_i \pp_{y_l}\Big\}f \\
&= P_1 (X_i f) - 2 P_1\Big( \sum_{l=1}^d (A_l x)_i \pp_{y_l} f\Big)\\
&= P_1(X_if) +2\sum_{l=1}^d P_1 \Big(\big\{(\tt  M^*)^{-1}(A_\cdot x)_i\big\}_l [X_{i_l}, X_{j_l}] f\Big).\end{split}\end{equation}
Next, for any $f\in C_0^\infty$, at point $(0,0)$ we have
\beg{equation*}\beg{split} &P_1\Big(\big\{(\tt M^*)^{-1}(A_\cdot x)_i\big\}_l [X_{i_l}, X_{j_l}] f\Big)\\
&= \int_{\R^{m+d}} p_1^0(x,y) \big\{(\tt M^*)^{-1}(A_\cdot x)_i\big\}_l (X_{i_l}X_{j_l}-X_{j_l}X_{i_l}) f(x,y) \d x \d y\\
&=  \int_{\R^{m+d}} p_1^0(x,y) (X_{i_l}f)(x,y) \Big[\big\{(\tt M^*)^{-1}(A_\cdot x)_i\big\}_l X_{j_l}\log p_1^0(x,y) +(\tt M^*)^{-1}(A_\cdot \si)_{ij_l}\Big]\d x\d y\\
&\quad -  \int_{\R^{m+d}} p_1^0(x,y) (X_{j_l}f)(x,y) \Big[\big\{(\tt M^*)^{-1}(A_\cdot x)_i\big\}_l X_{i_l}\log p_1^0(x,y) +\big\{(\tt M^*)^{-1}(A_\cdot \si)_{ii_l}\big\}_l\Big]\d x\d y\\
&= P_1\Big\{(X_{i_l}f)\Big[\big\{(\tt M^*)^{-1}(A_\cdot x)_i\big\}_l X_{j_l}\log p_1^0 +\big\{(\tt M^*)^{-1}(A_\cdot\si)_{ij_l}\big\}_l\Big]\Big\}\\
&\quad -P_1\Big\{(X_{j_l}f)\Big[\big\{(\tt M^*)^{-1}(A_\cdot x)_i\big\}_l X_{i_l}\log p_1^0 +\big\{(\tt M^*)^{-1}(A_\cdot\si)_{ii_l}\big\}_l\Big]\Big\}.\end{split}\end{equation*} Combining this with (\ref{PW}) and noting that
\beg{equation*}\beg{split} &  \sum_{l=1}^d \Big[\big\{(\tt M^*)^{-1}(A_\cdot x)_i\big\}_l X_{j_l}\log p_1^0 +\big\{(\tt M^*)^{-1}(A_\cdot\si)_{ij_l}\big\}_l\Big]^2\\
& \le
2\sum_{i,j=1}^m |(\tt M^*)^{-1}(A_\cdot \si)_{ij}|^2 + 4 \sup_{1\le l\le d} \big\{(\tt M^*)^{-1}(A_\cdot x)_i\big\}_l^2\GG(\log p_1^0)\end{split}\end{equation*} and the same holds for $i_l$ in place of $j_l$, we prove (\ref{GG'}) at point $(0,0).$

(3) Obviously, $C<\infty$ follows from
\beq\label{GG3} \int_{\R^{m+d}} |x|^2 p_1^0 (x,y) \GG(\log p_1^0) (x,y) \d x\d y<\infty.\end{equation}
Let $h\in C^\infty[0,\infty)$ such that $h|_{[0,1]}=1$ and $h|_{[2,\infty)}=0$. Let $f_n(x)= |x|^2 h(|x|/n)$. Then $f_n\in C_0^\infty(\R^{m+d}), n\ge 1.$ By Corollary \ref{C3.2} (4), there exists a constant $c>0$ such that
$$\GG(\log p_1^0)\le c\Big(1+ \ff{Lp_1^0}{p_1^0}\Big).$$ Combining these with $\big(L+\DDD+\ff{m+2}2 \big)p_1^0=0$ according to Lemma \ref{L4.2}, we arrive at
\beg{equation*}\beg{split} &\int_{\R^{m+d}} f_n(|x|) p_1^0 (x,y) \GG(\log p_1^0) (x,y) \d x\d y\\
&\le c \int_{\R^{m+d}} f_n(|x|) p_1^0 (x,y) \Big(1+ \ff{Lp_1^0}{p_1^0}\Big) (x,y) \d x\d y\\
&= c \int_{\R^{m+d}} f_n(|x|)Lp_1^0(x,y)\d x\d y + c P_1 |x|^2 (0,0) \\
&= cP_1|x|^2(0,0) -c\int_{\R^{m+d}} f_n(|x|) \Big(\DDD+\ff{m+2} 2\Big) p_1^0(x,y)\d x\d y\\
&= c P_1|x|^2(0,0)+ c P_1 (\DDD f_n)(0,0)\le c P_1|x|^2(0,0) + c P_1\Big(|x|^2 + \ff {\|h'\|_\infty} n |x|^3\Big)(0,0).\end{split} \end{equation*}
Letting $n\to\infty$ and noting that $P_1 |x|^p(0,0)<\infty $ holds for any $p\ge 1$, we obtain (\ref{GG3}).
\end{proof}

Finally,  we remark that extending the main result of \cite{Li} for the heat semigroup on the Heisenberg group (see also \cite[Section 5]{BBBC}), the following stronger estimate than (\ref{GG}) was proved in \cite{El} for the heat semigroup on a nilpotent Lie group of $H$-type:
 $$\ss{\GG(P_t f)}\le C P_t \ss{\GG(f)},\ \ f\in C_b^1(\R^3), t>0$$ for some constant $C>0$. This estimate implies the semigroup log-Sobolev inequality. However, in the moment we are not able to prove this type estimate under our more general framework. Note that to meet the requirement of $H$-type nilpotent Lie groups, in our framework  one has  to assume further that $m$ is even (see \cite[Proposition 2.1]{El})  and in \cite[(2.3)]{El} $J_{u_l}:= A_l$ is orthogonal.

\beg{thebibliography}{99}

 \bibitem{ATW06} M. Arnaudon, A. Thalmaier, F.-Y. Wang,
  \emph{Harnack inequality and heat kernel estimates
  on manifolds with curvature unbounded below,} Bull. Sci. Math. 130(2006), 223--233.


   \bibitem{BBBC} D. Bakry,  F. Baudoin,  M. Bonnefont, D. Chafa\"{\i}, \emph{ On gradient bounds for the heat kernel on the Heisenberg group,} J. Funct. Anal. 255 (                    2008),   1905--1938.

\bibitem{BB1} F. Baudoin, M. Bonnefont, \emph{Log-Sobolev inequalities for subelliptic operators satisfying a generalized curvature dimension inequality,} J. Funct. Anal. 262(2012), 2646--2676.

\bibitem{BB2} F. Baudoin, M. Bonnefont, N. Garofalo, \emph{A sub-Riemannian curvature-dimension inequality, volume doubling property and the Poincar\'e inequality,}
 arXiv:1007.1600.

\bibitem{BG} F. Baudoin, N. Garofalo, \emph{Curvature-dimension inequalities and Ricci lower bounds for sub-Riemannian manifolds with transverse symmetries,} arXiv:1101.3590.

\bibitem{BGM} F. Baudoin, M. Gordina, T. Melcher, \emph{Quasi-invariance for heat kernel measures on sub-Riemannian infinite-dimensional Heisenberg groups,} Trans. Amer. Math. Soc. 365(2013), 4313--4350.

\bibitem{B} J. M. Bismut, \emph{Large Deviations and the
Malliavin Calculus,} Boston: Birkh\"auser, MA, 1984.

\bibitem{CG} M. Cranston, A. Greven, \emph{Coupling and harmonic functions in the case of continuous time
Markov processes,} Stochastic Process Appl. 60(1995), 261--286.

\bibitem{CW} M. Cranston, F.-Y. Wang, F.-Y., \emph{A condition for the equivalence of coupling and shift-coupling,}
Ann. Probab. 28(2000), 1666--1679.



\bibitem{D} B. K. Driver,  \emph{Integration by parts for heat kernel measures revisited,} J. Math. Pures Appl. 76(1997), 703--737.

\bibitem{DG} B. K. Driver, M. Gordina, \emph{Integrated Harnack inequalities on Lie groups,} J. Diff. Geom. 83(2009), 501--550.

\bibitem{DM} B. K. Driver, T. Melcker, \emph{Hypoelliptic heat kernel inequalities on the Heisenberg group,} J. Funct. Anal. 221(2005), 340--365.

\bibitem{El} N. Eldredge, \emph{Gradient estimates for the subelliptic heat kernel on $H$-tpye groups,} 258(2010), 504--533.

\bibitem{GW} F.-Z. Gong, F.-Y. Wang, \emph{Heat kernel estimates with application to compactness of manifolds,} Quart. J. Math. 52
(2001), 171--180.

 \bibitem{Li} H.-Q. Li, \emph{Estimation optimale du gradient du semi-groupe de la chaleur sur le groupe de Heisenberg,} J. Funct. Anal. 236(2006), 369--394.
\bibitem{LY} P. Li, S. T. Yau, \emph{On the parabolic kernel of the Schr¡§odinger operator,} Acta Math. 156 (1986), 153-201.

\bibitem{W97} F.-Y. Wang,  \emph{Logarithmic Sobolev inequalities on noncompact Riemannian manifolds,} Probab. Theory Relat. Fields 109(1997), 417-424.

\bibitem{W10} F.-Y. Wang, \emph{Harnack inequalities on manifolds with boundary and applications,}    J.
Math. Pures Appl.   94(2010), 304--321.

\bibitem{WB} F.-Y. Wang, \emph{Analysis for Diffusion Processes on Riemannian Manifolds,} World Scientific, Singapore, 2013.

\bibitem{W12a} F.-Y. Wang, \emph{Derivative formula and gradient estimates for  Gruschin type semigroups,}    J. Theo. Probab. 27(2014), 80--95.

\bibitem{WN} F.-Y. Wang, \emph{Generalized curvature condition for subelliptic diffusion processes,} arXiv:1202.0778.

\end{thebibliography}
\end{document}